\documentclass[article, a4paper]{memoir}
\pdfoutput=1
\settrimmedsize{297mm}{210mm}{*}
\settypeblocksize{634pt}{448.13pt}{*}
\setulmargins{4cm}{*}{*}
\setlrmargins{*}{*}{1.5}
\setmarginnotes{17pt}{51pt}{\onelineskip}
\setheadfoot{\onelineskip}{2\onelineskip}
\setheaderspaces{*}{2\onelineskip}{*}
\checkandfixthelayout
\usepackage[hungarian]{babel}
\usepackage[T1]{fontenc}
\usepackage[utf8]{inputenc}
\PassOptionsToPackage{unicode}{hyperref}
\usepackage{hyperref, url, amssymb, amsmath, amsthm}
\hypersetup{pdfauthor={Gyula Magyarkuti},pdftitle={On the Hausdorff- and the Banach-Tarski paradox}}
\newtheorem{Tetel}{tétel}
\newtheorem{Allitas}[Tetel]{állítás}
\newtheorem{Definicio}[Tetel]{definíció}
\newtheorem{Lemma}[Tetel]{lemma}
\DeclareMathOperator{\id}{id}
\title{A Hausdorff- és a Banach\,--\,Tarski-paradoxonról}
\author{Magyarkuti Gyula\\
    Budapesti Corvinus Egyetem, Matematika Tanszék
}
\date{2020. December 17-én}
\begin{document}
\maketitle
\begin{abstract}
    An elementary approach to Banach-Tarski paradox is presented. 
    Very small amount of algebra and measure theory is required.
\end{abstract}
Egy narancs véges sok részre bontható olyan módon, hogy
a részekből két narancs is összerakható; vagy egy biliárd golyó és Stephan
Banach életnagyságú szobra feldarabolható páronként egybevágó részekre \cite{Stromberg-1979}. 
Ezek a szokásos népszerűsítő formái Hausdorff, Banach és Tarski idevágó tételeinek: 
\\
\emph{Adott két korlátos halmaz }$\mathbb{R}^{3}$\emph{-ban, amelyeknek a
belseje nem üres. Ekkor e két halmaz felbontható véges sok diszjunkt részhalmazra úgy, 
hogy az egyes részhalmazok egymásba vihetők a tér egybevágósági transzformációinak segítségével.}

A fenti állítás talán a 20. századi matematika legnagyobb meglepetése és jól
mutatja a háromdimenziós térbeli végtelen halmazok képzeletbeli jellegét. 
A tétel kizárja a nem triviális, végesen additív forgatás- és eltolás-invariáns 
mértékek létezését az $\mathbb{R}^{n}$ ($n\geq 3$) tér egész hatványhalmazán; 
és rámutat az olyan mérték kiterjesztési eljárások szükségszerűségére,
mint például a Lebesgue-mérték szokásos Caratheodory-féle konstrukciója.

\chapter*{Bevezetés}

Már az antikvitás kora óta ismert, hogy a végtelen fogalma nagyon gyorsan els%
ő látásra paradoxnak tűnő megállapításokhoz vezet: Bizonyos objektumok a mé%
retüket változtathatják látszólag mérettartó operációk hatására.

Galilei az 1638-ban megjelenő könyvében \cite{Galileo-1914} leírja azt a
megfigyelését, hogy pozitív egészeket és a négyzetszámokat kölcsönösen egymá%
shoz lehet rendelni, még akkor is ha a pozitív egészek négyzetszámokból és
nem négyzetszámokból állnak, így a pozitív egészek nyilván sokkal többen
vannak mint a négyzetszámok. Ebből levonja a következtetést: ``the
attributes `equal', `greater' and `less' are not applicable to infinite ...
quantities'', evvel előre sejtve majd 300 évvel későbbi mértékelméleti kutatá%
sok eredményeit bizonyos mértékek nem létezéséről.

Bár a végtelen részhalmazok paradoxnak érzett tulajdonságai a végtelen fogalm%
ának legelső használataikor is gondot okoztak, a paradox felbontások vizsgá%
lata csak a 20. század első éveiben a mértékelmélet formalizációjával kezdődö%
tt meg. Vitali által, a nem Lebesgue-mérhető halmazra adott klasszikus pé%
lda (1905), az első alkalom a paradox felbontások megjelenésére. Vitali
bizonyos mértékek létezésének képtelenségére akarta e példával a figyelmet
felhívni, és tíz évvel később Hausdorff a gömbfelszinen ,,konstruált'' való%
ban meglepő részhalmazokat, szintén bizonyos mértékek kizárásának igazolásá%
ul. Ez a munka inspirálta az 1920-as években Banach és Tarski kutatásait
akik egymástól függetlenül és különböző módszerrekkel igazolták lényegében a
fent kiemelt tételt.

Galilei megfigyelésének lényege, a megkettőzés. Kiindulva a pozitív egészekbő%
l két olyan halmazt készítünk, amelyek mindegyike ugyanakkora méretű mint a
kiindulásul vett halmaz. A megkettőzésnek ez az ötlete rejlik a jelen írás
sarokköveiben is (lásd: \ref{tolo}, \ref{2szabad}, \ref{btlemma}). Röviddel
azután, hogy Cantor számosság elmélete tisztázta a Galilei megfigyelésével
analog jelenségeket, Hausdorff, Banach és Tarski még bizarabb megkettőzések
lehetőségét fedezték fel, mivel azok pusztán a tér egybevágósági transzformá%
cióit használják. Ez az amit népszerűsítő írásokban úgy szoktak idézni, hogy
,,egy borsó szétszedhető véges sok darabra és összerakható eltolások és
elforgatások segítségével olyan nagyra, mint akár a nap''.

Tény, hogy a kiválasztási axiómát használjuk, ezért a fizikai megvalósítás
egy konkrét borsó esetében minimum kétségesnek tűnik. Wagon (1985) \cite
{Wagon-1985}-ben érdekes paradox konstruktív példákat találhatunk,
amelyekhez nem szükséges a kiválasztási axióma, viszont meg lehet mutatni,
hogy a kiválasztási axióma nélkül az állítás nem igaz. Mégis hibás úton jár,
aki a pusztán a kiválasztási axióma transzcendens mivoltát látja a jelenség
magyarázataként, mint ahogy Galilei megfigyelésének sincs köze a kiválasztá%
si axiómához. Például Dougherty és Foreman (1994) \cite
{DoughertyForeman-1994} eredménye szerint $\mathbb{R}^{n}$\emph{\ egységgömbjé%
ben található véges sok egymástól diszjunkt nyílt halmaz, amelyekből egybevág%
ósági transzformációkkal összerakható egy sűrű részhalmaza mondjuk a }$%
10^{100}$\emph{\ sugarú gömbnek}, és a kiválasztási axiómától ez az állítás f%
üggetlen!!! Sűrű halmaz persze mértékelméleti értelemben akár nagyon kicsi
is lehet, ami kétségbevonhatja a fenti állítás relevanciáját az eredeti probl%
émára nézve. Lásd még: \cite{Wikstrom-1995}

Logikusok is nekiestek a problámának és gyengítették a kiválasztási axiómát é%
ppen odáig, hogy B-T ne legyen, de a kiválasztási axióma valamennyi analí%
zisben használt ,,pozitív'' alkalmazására továbbra is lehetőség nyíljon. Ez
az úgynevezett \emph{Axiom of Dependent Choice}: 
\index{Axiom of Dependent Choice}\emph{Minden nem üres }$S$\emph{\ halmazhoz 
és tetszőleges }$f:S\rightarrow S$\emph{\ függvényhez létezik }$x_{n}\in S;$ 
$x_{n+1}=f\left( x_{n}\right) $\emph{\ tulajdonságú sorozat.} Schechter
(2000) \cite{Schechter-2000} megjegyzi, hogy evvel ekvivalens az ún. \emph{%
Dancs-Hegedűs-Medvegyev féle kiválasztási elv}. 
\index{Dancs-Hegedűs-Medvegyev féle kiválasztási elv}

Ilyen módon lehetőségünk van a Banach\,--\,Tarski-tételtől megszabadulni,
amennyiben az célunk. De legyen-e ez célunk? Ha igen, akkor felül kell vizsg%
álni más analízisbeli paradox jelenségeket is. Ekkor nem csak az a kérdés,
hogy hogyan fordulhatnak elő matematikai szörnyetegek, hanem az is kérdés,
hogy pontosan mik ezek. Folytonos függvény, amely sehol sem differenciálható
vajon kevésbé paradox jelenség? Evvel a gondolattal összefüggő filozófikus
hangvételű dolgozat Feferman (1998) \cite{Feferman-1998}.

Fontos tudni azt, amit itt nem igazolok, hogy $\mathbb{R}$-ben és $\mathbb{R}^{2}$%
-ben van \emph{Banach-mérték}, 
\index{Banch-mérték}tehát olyan $\mathcal{P}\left( \mathbb{R}^{j}\right)
\rightarrow \mathbb{R}$ eltolás-, illetve forgatás- invariáns végesen additív
halmazfüggvény ($j=1,2$), amely az egész hatványhalmazon van értelmezve és
az egységgömb mértéke 1. Világos, hogy ekkor B-T nem lehet, tehát B-T és
Banach-mérték létezése egymással alternatív állítások. Lásd: \cite
{Wagon-1985}

Jelen írás célja összesen annyi, hogy a bevezető előtt kiemelt tételt, a
Banach\,--\,Tarski-tétel erős alakját bebizonyítsuk. Csak elemi eszközöket
fogunk használni. Az általunk használt legerősebb eszköz a három dimenziós té%
r ortogonális transzformációinak karakterizációja. Az állítást az áttekinthet%
őség kedvéért konkrétan $\mathbb{R}^{3}$-ban igazolom, de semmi különöset nem
kellene hozzátenni ahhoz, hogy $\mathbb{R}^{n}$, $n>3$ mellett is belássuk.

\section*{Jelölések}

Jelölje az egész írásban $B^{n+1}$ az $\mathbb{R}^{n+1}$ tér euklideszi-normá%
val számított egységgömbjét és $S^{n}$ annak felületét. Legyen $SO\left(
n+1\right) $ az $S^{n}$ felület forgatás csoportja, azaz az $\mathbb{R}^{n}$ té%
r azon $A$ ortogonális ($A^{\star }=A^{-1}$) transzformációinak $S^{n}$-re
való leszűkítése, amelyekre $\det A=1.$ Jelölje $G_{n}$ az $\mathbb{R}^{n}$
izometria csoportját. Tudjuk, hogy $G_{n}$ eltolások és forgatások kompozíció%
iból áll.

Általában $G$-vel multiplikatív irásmódú csoportot fogunk jelölni, és ekkor $%
1$-el ennek reprodukáló elemét jelöljük.

A $\sigma $ és $\tau $ két elem generálta szabad csoporton a $\sigma ,\sigma
^{-1},\tau ,\tau ^{-1}$ négy különböző jelből előállítható redukált szavak
halmazát értjük, ahol a csoport művelet a konkatenáció (egymás mögé írás).
Az identikus elemet, tehát a $0$ hosszú redukált szót $%
\id$ fogja jelölni. Redukált szón egy $\alpha _{1}^{i_{1}}\alpha
_{2}^{i_{2}}\ldots \alpha _{n}^{i_{n}}$ alakú kifejezést értünk, ahol $%
i_{1},\ldots ,i_{n}$ pozitív egészek, $\alpha _{i}\in \left\{ \sigma ,\sigma
^{-1},\tau ,\tau ^{-1}\right\} $ és $i<n$ esetén $\alpha _{i}\neq \alpha
_{i+1}^{-1}.$ Perszer e szabad csoporttal izomorf csoportot is szabad
csoportnak fogunk mondani.

A szokásoknak megfelelően $X$ egyenlőre egy absztrakt halmaz, és $\mathcal{P}%
\left( X\right) $ annak hatványhalmaza.

\chapter{A Hausdorff-paradoxon}

Először Vitali programját gondoljuk át, persze sokkal absztraktabb kö%
rnyezetben, majd bevezetjük a kongruens, az ekvidekompozábilis és a paradox
halmazok és csoportok fogalmát. Pédát is adunk paradox halmazra, amiből már k%
önnyen fog következni, hogy $S^{2}$ egy megszámlálható halmaztól eltekintve
megduplázható.

\section{Particionálás egy csoport hatásaként}

\begin{Definicio}[a $G$ csoport hat az $X$ halmazon]
\index{csoport hatása egy halmazon} Legyen $G$ egy csoport és $X$ egy rögzí%
tett halmaz. Azt mondjuk, hogy $G$\emph{\ az }$X$\emph{-en hat}, ha minden $%
g\in G$ -hez létezik egyértelműen egy $%
\widehat{g}:X\rightarrow X$ bijekció, melyre 
\begin{equation}
\widehat{gh}=\widehat{g}\circ \widehat{h}\text{\quad és\quad }\widehat{1}=%
\id.  \label{hat}
\end{equation}

Azt mondjuk, hogy $G$ \emph{fixpontmentesen hat} az $X$-en, ha minden $g\in
G,$ $g\neq 1$ esetén $\widehat{g}$ egy fixpontmentes bijekciója $X$-nek. 
\index{csoport fixpontmentes hatása egy halmazon}
\end{Definicio}

Amennyiben $G$ az $X$-en hat, úgy a 
\[
\widehat{g^{-1}}=\left( \widehat{g}\right) ^{-1}
\]
azonosság is teljesül, hiszen $g\cdot g^{-1}=1$ miatt $\widehat{g}\circ 
\widehat{g^{-1}}=\id$ is fennáll.

\label{csillagtrafo}Legyen $G$ egy az $X$ halmazon ható csoport, és $%
M\subseteq X$ egy rögzített részhalmaz. Vezessük be, a következő $\mathcal{P}%
\left( G\right) \rightarrow \mathcal{P}\left( X\right) $ operációt: 
\[
B\subseteq G\text{ esetén }B^{\star }\doteq \left\{ x\in X:\exists m\in M%
\text{, }\exists g\in B,x=\hat{g}\left( m\right) \right\} =\cup \left\{ 
\hat{g}\left( M\right) :\hat{g}\in B\right\} .
\]
Jelölje tetszőleges $B\subseteq G$ és tetszőleges $g\in B$ mellett $gB\doteq
\left\{ gh:h\in B\right\} \subseteq X$ halmazt. Ekkor 
\[
\left( gB\right) ^{\star }=\hat{g}\left( B^{\star }\right) ,
\]
ugyanis 
\[
\hat{g}\left( B^{\star }\right) =\cup \left\{ \hat{g}\left( \hat{h}\left(
M\right) \right) :h\in B\right\} =\cup \left\{ \widehat{gh}\left( M\right)
:h\in B\right\} \cup \left\{ \widehat{f}\left( M\right) :f\in gB\right\}
=\left( gB\right) ^{\star }.
\]

\label{csoport}Érdemes a következő példákat a későbbi tételek folyamán is
sze melőtt tartani.

\begin{enumerate}
\item  Legyen $G$ egy csoport és $X\doteq G.$ Adott $g\in G$ mellett a $%
\widehat{g}:G\rightarrow G$ bijekció legyen a $g$-vel képzett balszorzás,
azaz 
\[
\widehat{g}\left( h\right) \doteq gh.
\]
Könnyen látható, hogy ez valóban bijekció és (\ref{hat}) is fennáll, tehát $G
$ valóban egy a $G$-n ható csoport. Sőt: $G$ fixpontmentesen hat $G$-n.

\item  Legyen $n\geq 2$ természetes szám mellett $G=SO\left( n+1\right) $,
azaz az $S^{n}$ forgatásainak csoportja és $X=S^{n}.$ Adott $g\in G$ esetén $%
\widehat{g}=g.$ Látható, hogy $SO\left( n+1\right) $ hat az $S^{n}$-en, de
nem fixpontmentesen. Ugyanis $n\geq 2$ miatt, a forgatás tengelye és az $n+1$%
-dimenziós egységgömb két közös pontja egyben a forgatás fixpontja is.

\item  Legyen $G=SO\left( 2\right) $, azaz az $\mathbb{R}^{2}$ sík $S^{1}$ egysé%
gkörének forgatáscsoportja és $X=S^{1}.$ Adott $g\in G$ esetén $\widehat{g}%
=g.$ Az $SO\left( 2\right) $ csoport fixpontmentesen hat az $S^{1}$-en.

\item  Legyen $G$ az $\mathbb{R}^{2}$ sík $S^{1}$ egységkörének racionális szö%
gekkel való forgatáscsoportja, és $X=S^{1}.$ Adott $g\in G$ esetén $\widehat{%
g}\doteq g.$ A $G$ csoport persze fixpontmentesen hat az $S^{1}$-en.

\item  Legyen $X=\left[ 0,2\pi \right) ,$ és $G$ a $\left[ 0,2\pi \right) $%
-beli racionális számok additív csoportja a moduló $2\pi $-vel számított ö%
sszeadás művelettel ellátva. Adott $g\in G$ esetén $\widehat{g}:\left[
0,2\pi \right) \rightarrow \left[ 0,2\pi \right) $ a 
\[
\widehat{g}\left( r\right) \doteq g+r\text{\quad }\mod\left( 2\pi
\right) 
\]
módon megadott függvény. Világos, hogy ez csak átfogalmazása az előző példá%
nak, és $G$ fixpontmentesen hat $\left[ 0,2\pi \right) $-en.

\item  Legyen $X=\mathbb{R}$ és $G$ a racionális számok additív csoportja.
Adott $g\in G$ esetén $\widehat{g}:\mathbb{R}\rightarrow \mathbb{R}$ a 
\[
\widehat{g}\left( r\right) \doteq g+r
\]
módon definiált függvény. Világos, hogy $G$ fixpontmentesen hat $X$-en.
\end{enumerate}

\begin{Allitas}
\label{particioallitas}Legyen $G$ egy az $X$ halmaz fölött ható csoport. Jelö%
lje $G\left( x\right) \doteq \left\{ \widehat{g}\left( x\right) :g\in
G\right\} $ az $x$ pont pályáját.

$\left( 1\right) $ Ekkor a pályákból álló halmazrendszer az $X$ alaphalmaz
egy partícióját alkotja, azaz 
\begin{equation}
\left\{ G\left( x\right) :x\in X\right\} \text{ partíciója }X\text{-nek.}
\label{particio}
\end{equation}

$\left( 2\right) $ Tegyük most fel, hogy $G$ fixpontmentesen hat az $X$
halmazon. Vegyünk ki a fenti partíció minden halmazából pontosan egy-egy
elemet és jelölje $M$ az így kapott reprezentáns halmazt. Ekkor a $\left\{ 
\widehat{g}\left( M\right) :g\in G\right\} $ halmazrendszer partíciója $X$%
-nek, amely ekvipotens $G$-vel, azaz 
\begin{equation}
\hat{g}\left( M\right) \neq \varnothing ;\quad g_{1}\neq g_{2}\Rightarrow 
\hat{g}_{1}\left( M\right) \cap \hat{g}_{2}\left( M\right) =\varnothing
;\quad X=\cup _{g\in G}\hat{g}\left( M\right) .  \label{particio2}
\end{equation}

$\left( 3\right) $ Rögzítsük az iménti $M$ reprezentáns halmazt, és tekintsü%
nk a $G$ egy $G=\cup _{i=}^{n}B_{i}$ partícióját. Ekkor $X=\cup
_{i=1}^{n}B_{i}^{\star }$ az $X$ partícióját adja ahol $\star $ a \ref
{csillagtrafo}. szakaszban bevezetett transzformáció.
\end{Allitas}

\begin{proof}
$\left( 1\right) $ Először is $x\in G\left( x\right) $ minden $x\in X$
mellett fennáll, hiszen az $1\in G$ multiplikatív egységelemre $\widehat{1}=%
\id.$ Másodszor, ha $u\in G\left( x\right) \cap G\left( y\right) $
azaz valamely $r,s\in G$ mellett $\widehat{r}\left( x\right) =u=\widehat{s}%
\left( y\right) $ és $z\in G\left( x\right) ,$ akkor $z\in G\left( y\right) $
is fennáll. Ugyanis ha $z=\widehat{t}\left( x\right) ,$ valamely $t\in G$
esetén, akkor 
\[
z=\widehat{t}\left( x\right) =\widehat{t}\left( \widehat{r}^{-1}\left(
u\right) \right) =\widehat{t}\left( \widehat{r^{-1}}\left( u\right) \right) =%
\widehat{t}\left( \widehat{r^{-1}}\left( \widehat{s}\left( y\right) \right)
\right) =\widehat{\left( tr^{-1}s\right) }\left( y\right) \in G\left(
y\right) .
\]
Megmutattuk tehát, hogy $G\left( x\right) \subseteq G\left( y\right) .$
Hasonlóan kapjuk, hogy $G\left( y\right) \subseteq G\left( x\right) ,$ tehát 
\[
G\left( x\right) \cap G\left( y\right) \neq \varnothing \Rightarrow G\left(
x\right) =G\left( y\right) ,
\]
ami azt jelenti, hogy a (\ref{particio}) valóban egy partícióját alkotja $X$%
-nek.

\noindent $\left( 2\right) $ Most tegyük fel, hogy $\left\{ \widehat{g}:g\in
G,g\neq 1\right\} $ fixpontmentes bijekciók. Persze minden $x\in X$ mellett $%
x\in G\left( x\right) $ és arra az egyetlen $m\in M$-re, amelyre $m\in
G\left( x\right) $, az $x\in G\left( x\right) =G\left( m\right) $ is teljesü%
l, ami épp azt jelenti, hogy létezik $g\in G$, melyre $x=\widehat{g}\left(
m\right) $.\newline
Tegyük fel most, hogy valamely $g,h\in G$ mellett $\widehat{g}\left(
M\right) \cap \widehat{h}\left( M\right) \neq \varnothing $, azaz létezik $%
m,m^{\prime }\in M,$ melyekre $\widehat{g}\left( m\right) =\widehat{h}\left(
m^{\prime }\right) .$ Ekkor persze $G\left( m\right) \cap G\left( m^{\prime
}\right) \neq \varnothing ,$ ezért $G\left( m\right) =G\left( m^{\prime
}\right) $. De mivel a $\left( \text{\ref{particio}}\right) $ partíció
minden részhalmazából csak egyetlen elemet vettünk ki az $M$ halmaz képzésé%
hez, ezért $m=m^{\prime }$ is teljesül. Így $\widehat{g}\left( m\right) =%
\widehat{h}\left( m\right) ,$ amiből $m=\left( \widehat{h}\right)
^{-1}\left( \widehat{g}\left( m\right) \right) =\widehat{\left(
h^{-1}g\right) }\left( m\right) $ következik. Tehát a $h^{-1}g\in G$ egy
olyan elem, melyhez rendelt bijekciónak van fixpontja. A feltételünk szerint
ez csak úgy lehetséges, ha $h^{-1}g=1,$ azaz $h=g.$

\noindent $\left( 3\right) $ Az állítás szinte csak átfogalmazása $\left(
2\right) $-nek.
\end{proof}

A bizonyítandó állítás szempontjából kitérő, de érdemes felfigyelni arra,
hogy a fenti (\ref{particioallitas}) állítás szoros kapcsolatban van a nem mé%
rhető halmaz szokásos konstrukciójával.

\index{nem mérhető halmaz konstrukciók} Tekintsük először a $4.$ vagy $5.$ pé%
ldát. Az ott konstruált $M$ halmaz biztosan nem Lebesgue-mérhető, hiszen az
előző állítás szerint a $\left[ 0,2\pi \right) $ intervallum előáll 
\[
\left[ 0,2\pi \right) =\bigcup_{g\in \mathbb{Q}}\left( g+M\right) 
\]
diszjunkt egyesítés alakban. Mivel a Lebesgue-mérték eltolás invariáns, ezé%
rt $M$ mérhetőségéből a $g+M$ eltolt halmazok mérhetősége is következne,
majd a $2\pi =\lim_{n\rightarrow \infty }n\lambda \left( M\right) ,$ ami
nyilván ellentmondás, hiszen ez utóbbi szám csak $0$ vagy $+\infty $ lehet.

A $6.$ példa alkalmazása is $\mathbb{R}$ egy nem Lebesgue-mérhető halmazához
vezet. A (\ref{particioallitas}) állításbeli $M$ halmaz konstrukcióját csak
egy kicsit kell megváltoztatnunk. Egyrészt minden $x\in \mathbb{R}$ mellett van
olyan $y\in \mathbb{R}$, $0\leq y\leq 1,$ hogy $x-y\in \mathbb{Q.}$ Ez azt
jelenti, hogy $G\left( x\right) \cap \left[ 0,1\right] \neq \varnothing $
minden $x\in \mathbb{R}$ esetén, tehát a reprezentánsokat választhatjuk a $%
\left[ 0,1\right] $ intervallumból, így az $M\subseteq \left[ 0,1\right] $
tartalmazás feltehető. Másrészt 
\[
\left[ 0,1\right] \subseteq \bigcup_{g\in \mathbb{Q}  \left| g\right|
\leq 2}\left( g+M\right) \doteq H,
\]
hiszen $x\in \left[ 0,1\right] $-hez létezik $g\in \mathbb{Q}$ melyre $x\in g+M,
$ azaz $x-g=m\in \left[ 0,1\right] .$ Ezért $\left| g\right| =\left|
x-m\right| \leq 2.$ Azt kaptuk tehát, hogy a $H$ halmazra $\left[ 0,1\right]
\subseteq H\subseteq \left[ -2,2\right] +\left[ 0,1\right] =\left[ -2,3%
\right] .$ Amennyiben tehát $M$ Lebesgue-mérhető lenne, akkor $H$ is
Lebesgue-mérhető lenne, és $H$ Lebesue-mértékére $1\leq \lambda \left(
H\right) \leq 5$ állna fenn. De $H$ definíciója szerint csak $\lambda \left(
H\right) =0,$ vagy $\lambda \left( H\right) =\infty $ lehetséges, ami ugyan
azt az ellentmondást adja, mint az előző példa.

A fenti példából az is világos, hogy nincs $\mathcal{P}\left( S^{1}\right) $%
-en $\mu $ megszámlálhatóan additív forgatás-invariáns mérték, amelyre $%
0<\mu \left( S^{1}\right) <\infty .$ Ebből azonnal következi, hogy nincs
olyan valós, eltolás invariáns $\mu $ mérték az egész $\mathcal{P}\left( 
\mathbb{R}\right) $-en, amelyre $0<\mu \left( \left[ 0,1\right] \right) <1,$ ezé%
rt az $A\mapsto A\times \left[ 0,1\right] ^{n-1}$ leképezés segítségével rögt%
ön láthatjuk, hogy nincs eltolás invariáns $\mu $ mérték a $\mathcal{P}%
\left( \mathbb{R}^{n}\right) $ halmazon értelmezve, amelyre $\mu \left( \left[
0,1\right] ^{n}\right) =1.$

A (\ref{particioallitas}) állításra úgy is nézhetünk, mint egy általános eszk%
öz nem mérhető halmazok konstrukciójához.

\section{Az ekvidekompozábilis halmazok}

\begin{Definicio}[kongruens halmazok, $A\cong _{G}B$]
\index{kongruens halmazok} Legyen $G$ egy az $X$ halmazon ható csoport,
valamint $A,B\subseteq X.$ Azt mondjuk, hogy az $A$\emph{\ halmaz }$G$\emph{%
-kongruens a }$B$\emph{\ halmazzal} ($A\cong _{G}B$), ha létezik $g\in G,$
amelyre $%
\widehat{g}\left( A\right) =B.$
\end{Definicio}

\begin{Definicio}[ekvidekompozábilis halmazok, $A\sim _{G}B$]
\index{ekvidekompozábilis halmazok} Legyen $G$ egy olyan csoport, amely az $X
$ halmazon hat, továbbá $A,B\subseteq X.$ Azt mondjuk, hogy az $A$\emph{\ és
a }$B$\emph{\ halmazok }$G$\emph{-ekvidekompozábilisek} ($A\sim _{G}B$), ha
létezik $A$-nak és $B$-nek olyan véges partíciója, amelynek elemei páronként 
$G$-kongruensek egymással.

Magyarul: létezik $A=\cup _{i=1}^{n}A_{i}$ és $B=\cup _{i=1}^{n}B_{i}$
diszjunkt felbontás, amelyekre $A_{i}\cong _{G}B_{i}$ minden $i=1,\ldots ,n$
mellett.
\end{Definicio}

Elképzelhető, hogy egy halmaz ekvidekompozábilis egy valódi részhalmazával:

\begin{Allitas}
\label{tolo}Legyen $G$ egy az $X$ halmazon ható csoport, valamint $%
P\subseteq A\subseteq X,$ és tegyük fel, hogy létezik $g\in G$ és létezik $%
Q\subseteq X,$ amelyekre 
\[
P\subseteq Q\subseteq A\quad 
\text{valamint\quad }\hat{g}\left( Q\right) =Q\smallsetminus P.
\]
Ekkor 
\[
A\sim _{G}\left( A\smallsetminus P\right) .
\]
\end{Allitas}

\begin{proof}
Világos, hogy az $A=Q\cup \left( A\smallsetminus Q\right) $ valamint az $%
A\smallsetminus P=\left( Q\smallsetminus P\right) \cup \left(
A\smallsetminus Q\right) $ diszjunkt egyesítések, valamint $\hat{g}\left(
Q\right) =Q\smallsetminus P$ és $\hat{1}\left( A\smallsetminus Q\right)
=A\smallsetminus Q.$
\end{proof}

A fenti állítás feltételeit kielégítő konkrét transzformációkat találunk pé%
ldaként a (\ref{btlemma}) lemmában és a (\ref{2szabad}) állításban.

\begin{Definicio}[paradox halmaz]
\index{paradox halmaz} Legyen $G$ egy csoport, amely az $X$ halmazon hat és $%
A\subseteq X$. Azt mondjuk, hogy az $A$\emph{\ halmaz }$G$\emph{-paradox},
ha létezik léteznek $B,C\subseteq A;$ $B\cap C=\varnothing $ diszjunkt ré%
szhalmazok, amelyek külön-külön $G$-ekvidekompozábilisek $A$-val, azaz $%
B\sim _{G}A$ és $C\sim _{G}A$ is fennáll.
\end{Definicio}

Az (\ref{paradox}) állításban látjuk majd, hogy a fenti definícióban $B\cup
C=X$ is feltehető. Ehhez a Banach\,--\,Schröder\,--\,Bernstein-tétel kell, de amig
nincs szükségünk rá addig nem akarom kihasználni.

\begin{Definicio}[paradox csoport]
\index{paradox csoport} Egy csoportot \emph{paradoxnak }mondunk, ha úgy
tekintve mint egy saját maga felett ható csoport (lásd a \ref{csoport}.
szakasz 1. példát), a $G$ halmaz $G$-paradox.
\end{Definicio}

Az alábbiakban példát adunk paradox csoportra. Ha arra gondolunk, hogy az egé%
sz számok azért ekvipotensek a páros számokkal mert mindegyiket meg kell
szoroznunk egy $2$-sel, talán nem is olyan meglepő, hogy milyen könnyen talá%
lkozhatunk egy ártatlan paradox csoporttal.

\begin{Allitas}
\label{2szabad}%
\index{szabad csoport} Minden két elem generálta szabad csoport paradox.
\end{Allitas}

\begin{proof}
Szabad csoport elemeinek redukált szóként való előállítása egyértelmű, é%
rtelmes tehát a következő halmaz megadása. Jelölje $W\left( \alpha \right) $
a $G$ csoport $\alpha $-val kezdődő redukált szavaiból álló részhalmazát,
amennyiben $\alpha =\sigma ^{\pm 1}$ vagy $\alpha =\tau ^{\pm 1},$ ahol $%
\sigma $ és $\tau $ a két generáló eleme a $G$ szabad csoportnak. Világos,
hogy $G$ előáll az alábbi diszjunkt egyesítésként: 
\[
G=W\left( \sigma \right) \cup W\left( \sigma ^{-1}\right) \cup W\left( \tau
\right) \cup W\left( \tau ^{-1}\right) \cup \left\{ 
\id\right\} .
\]
Jelölje $P\doteq W\left( \tau \right) \cup W\left( \tau ^{-1}\right) \cup
\left\{ \id\right\} $ és $Q\doteq P\cup W\left( \sigma \right) .$ Lá%
tható, hogy $\sigma \left( Q\right) =W\left( \sigma \right) =Q\smallsetminus
P,$ ami (\ref{tolo}) állítás szerint azt jelenti, hogy $G\sim W\left( \sigma
\right) \cup W\left( \sigma ^{-1}\right) .$ A fenti gondolatban $\sigma $
helyett $\tau $-t írva, azt kapjuk, hogy $G\sim W\left( \tau \right) \cup
W\left( \tau ^{-1}\right) ,$ ami e két utóbbi halmaz diszjunktsága miatt
igazolja az állítást.
\end{proof}

\begin{Allitas}
\label{partranz}Tekintsünk egy $A\subseteq X$ halmazt, amely $G$-paradox és
legyen $A^{\prime }\subseteq X;$ $A\sim _{G}A^{\prime }.$ Ekkor $A^{\prime }$
is $G$-paradox halmaz.
\end{Allitas}

\begin{proof}
Tegyük fel, hogy $B,C\subseteq A$ diszjunkt halmazok, amelyekre $B\sim _{G}A$
és $C\sim _{G}A$ is fennáll. Létezik $A=\cup _{k=1}^{n}A_{k}$ partíció és lé%
teznek $g_{k}\in G$ elemek amelyekkel $A^{\prime }=\cup _{k=1}^{n}\hat{g}%
_{k}\left( A_{k}\right) $ az $A^{\prime }$ egy partícióját alkotja. Legyen $%
B^{\prime }\doteq \cup _{k=1}^{n}g_{k}\left( A_{k}\cap B\right) $ és hasonló%
an $C^{\prime }\doteq \cup _{k=1}^{n}g_{k}\left( A_{k}\cap C\right) .$ Vilá%
gos, hogy $B^{\prime }\cap C^{\prime }=\varnothing $ és $B\sim _{G}B^{\prime
}$ valamint $C\sim _{G}C^{\prime }$. A $\sim _{G}$ reláció tranzitivitása
miatt tehát $B^{\prime }\sim _{G}A^{\prime }$ és $C^{\prime }\sim
_{G}A^{\prime }.$ Ezt kellett belátni.
\end{proof}

\section{A Hausdorff-paradoxon}

Először is azt kell látnunk, hogy amennyiben van egy paradox csoportunk
akkor ennek segítségével könnyen kaphatunk más paradox halmazokat is. Sőt,
az alábbi állítás megfordítása is igaz (lásd \cite{Su-1990} és \cite
{Wagon-1985}), de itt most érdektelen.

\begin{Allitas}
\label{hausdorff0}Ha $G$ az $X$ halmazon fixpontmentesen ható paradox
csoport, akkor $X$ is $G$-paradox.
\end{Allitas}

\begin{proof}
Legyenek $B,C\subseteq G$ diszjunkt részhalmazai $G$-nek, amelyek külön-külö%
n ekvidekompozábilisek $G$-vel. Léteznek tehát 
\[
B=\cup _{i=1}^{n}B_{i}\quad \ \text{és\quad }\ C=\cup _{j=1}^{m}C_{j}
\]
partíciók, valamint $g_{i},h_{j}\in G$ elemek $i=1,\ldots ,n$, és $%
j=1,\ldots ,m$, amelyekkel 
\[
\cup _{i=1}^{n}g_{i}B_{i}=G\quad \ \text{és\quad }\ \cup _{j=1}^{m}h_{j}Cj=G
\]
is egy-egy partícióját alkotja $G$-nek. A (\ref{particioallitas}) állításból
tudjuk, hogy ha $M$ a $\left\{ G\left( x\right) :x\in X\right\} $ partíció
tetszőleges reprezentánsaiból álló halmaz, akkor $\left\{ \widehat{g}\left(
M\right) :g\in G\right\} $ halmazrendszer is partíciója $X$-nek, ezért a $G$
egy particionálása az $X$ egy particionálását indukálja a (\ref{csillagtrafo}%
) szakaszban bevezetett $\star :\mathcal{P}\left( G\right) \rightarrow 
\mathcal{P}\left( X\right) $ transzformáció segítségével. Tehát $B^{\star
}=\cup B_{i}^{\star }$ particionálását adja $B^{\star }$-nak és $G^{\star
}=\cup \left( g_{i}B_{i}\right) ^{\star }$ particionálását adja $G^{\star }$%
-nak. De 
\[
B^{\star }=\cup _{i=1}^{n}B_{i}^{\star }\quad \sim _{G}\quad \cup _{i=1}^{n}%
\hat{g}_{i}\left( B_{i}^{\star }\right) =\cup _{i=}^{n}\left(
g_{i}B_{i}\right) ^{\star }=\left( \cup g_{i}B_{i}\right) ^{\star }=G^{\star
}=X.
\]
Hasonlóan $C^{\star }\sim _{G}X$ is fennáll.
\end{proof}

A probléma tehát abban áll, hogy milyen módon kapunk $SO\left( 3\right) $
forgatás csoport két elem generálta szabad részcsoportját?

\begin{Allitas}[Świerczkowski (1958)]
\index{forgatás csoport} Az $S^{2}$ gömbfelelület $SO\left( 3\right) $ forgat%
ás csoportjának van két elem generálta szabad részcsoportja.
\end{Allitas}

\begin{proof}
Azt kell megmutatnunk, hogy van $SO\left( 3\right) $-nak két eleme -- $\phi 
$ és $\rho $ --, amelyek által generált részcsoportban a redukált szavak elő%
állítása egyértelmű. Ez avval ekvivalens, hogy a $\phi ^{\pm 1}$ és a $\rho
^{\pm 1}$ jelekből az identitás transzformáció legalább egy hosszú szóként
nem rakható ki.

Legyen $\phi $ a $z$-tengely körüli $\arccos 
\frac{1}{3}$ fokkal való forgatás és $\rho $ az $x$-tengely körüli
ugyanekkora szögű forgatás. Ekkor $\phi $-nek és $\phi ^{-1}$-nek valamint $%
\rho $-nak és $\rho ^{-1}$-nek a mátrix reprezentációja: 
\[
\phi ^{\pm 1}=\left( 
\begin{array}{ccc}
\frac{1}{3} & \mp \frac{2\sqrt{2}}{3} & 0 \\ 
\pm \frac{2\sqrt{2}}{3} & \frac{1}{3} & 0 \\ 
0 & 0 & 1
\end{array}
\right) \text{\quad }\rho ^{\pm 1}=\left( 
\begin{array}{ccc}
1 & 0 & 0 \\ 
0 & \frac{1}{3} & \mp \frac{2\sqrt{2}}{3} \\ 
0 & \pm \frac{2\sqrt{2}}{3} & \frac{1}{3}
\end{array}
\right) .
\]
Tekintsünk tehát a $\phi $ és $\rho $ generálta részcsoportban egy $w$ reduká%
lt szót, amely egyenlő az identitás transzformációval és legalább egy hosszú%
. Feltehető, hogy $w$ utolsó betűje $\phi ,$ hiszen ezt elérhetjük a $\phi
^{-1}w\phi $ konjugálással. \textit{Megmutatjuk, hogy amennyiben }$w$\textit{%
\ egy a }$\phi ^{\pm 1}$\textit{\ és a }$\rho ^{\pm 1}$\textit{\ elemekből
alkotott redukált redukált szó, melynek utolsó betűje }$\phi $\textit{,
akkor } 
\[
w\left[ 1,0,0\right] =\left[ a,b\sqrt{2},c\right] /3^{k}
\]
\textit{alakú, ahol }$a,b,c\in \mathbb{Z}$\textit{\ egészek, továbbá }$3$%
\textit{\ nem osztója }$b$\textit{-nek.} Ebből bőven következik, hogy $w$
nem lehet az identitás transzformáció. (Látjuk majd, hogy a gyengébb állítá%
st nehezebb lenne indukcióval bizonyítani.)

A $w$ hossza szerinti indukció következik. Amennyiben $w$ egy hosszú, azaz
csak $\phi $-ből áll, akkor 
\[
w\left[ 1,0,0\right] =\phi \left[ 1,0,0\right] =\left[ 1,2\sqrt{2},0\right]
/3 
\]
ami valóban megfelel az állításnak, 
\[
a_{1}\doteq 1;\quad b_{1}\doteq 2;\quad c_{1}\doteq 0 
\]
jelöléssel.

Tegyük fel most, hogy $k>1$ és minden $k$-nál rövidebb redukált szó teljesí%
ti állításunkat, valamint $w$ éppen $k$ hosszú. Ekkor 
\[
w=\phi ^{\pm 1}w^{\prime }\text{\quad vagy\quad }w=\rho ^{\pm 1}w^{\prime }, 
\]
ahol $w^{\prime }\left[ 1,0,0\right] =\left[ a_{k-1},b_{k-1}\sqrt{2},c_{k-1}%
\right] /3^{k-1}$ alakban írható fel, persze $a_{k-1},b_{k-1},c_{k-1}\in 
\mathbb{Z}$ és $b_{k-1}$ nem osztható $3$-al. Amennyiben $w=\phi ^{\pm
1}w^{\prime }$ alakú, akkor az alábbi egyszerű számolgatás szerint: 
\begin{eqnarray*}
3^{k-1}w\left[ 1,0,0\right] &=&3^{k-1}\phi ^{\pm 1}w^{\prime }\left[ 1,0,0%
\right] =\phi ^{\pm 1}\left[ a_{k-1},b_{k-1}\sqrt{2},c_{k-1}\right] \\
&=&a_{k-1}\left[ \frac{1}{3},\pm \frac{2\sqrt{2}}{3},0\right] +b_{k-1}\sqrt{2%
}\left[ \mp \frac{2\sqrt{2}}{3},\frac{1}{3},0\right] +c_{k-1}\left[ 0,0,1%
\right] \\
&=&\left[ a_{k-1}\mp 4b_{k-1},\pm 2\sqrt{2}a_{k-1}+\sqrt{2}b_{k-1},3c_{k-1}%
\right] /3=\left[ a_{k},b_{k}\sqrt{2},c_{k}\right] /3,
\end{eqnarray*}
ahol tehát 
\begin{equation}
a_{k}\doteq a_{k-1}\mp 4b_{k-1};\quad b_{k}\doteq b_{k-1}\pm 2a_{k-1};\quad
c_{k}=3c_{k-1}.  \tag{$\dag $}
\end{equation}
Hasonlóan, ha $w=\rho ^{\pm }w^{\prime }$, akkor 
\begin{eqnarray*}
3^{k-1}w\left[ 1,0,0\right] &=&3^{k-1}\rho ^{\pm 1}w^{\prime }\left[ 1,0,0%
\right] =\rho ^{\pm 1}\left[ a_{k-1},b_{k-1}\sqrt{2},c_{k-1}\right] \\
&=&a_{k-1}\left[ 1,0,0\right] +b_{k-1}\sqrt{2}\left[ 0,\frac{1}{3},\pm \frac{%
2\sqrt{2}}{3}\right] +c_{k-1}\left[ 0,\mp \frac{2\sqrt{2}}{3},\frac{1}{3}%
\right] \\
&=&\left[ 3a_{k-1},\sqrt{2}b_{k-1}\mp 2\sqrt{2}c_{k-1},\pm 4b_{k-1}+c_{k-1}%
\right] /3=\left[ a_{k},b_{k}\sqrt{2},c_{k}\right] /3,
\end{eqnarray*}
ahol most 
\begin{equation}
a_{k}\doteq 3a_{k-1};\quad b_{k}\doteq b_{k-1}\mp 2c_{k-1};\quad c_{k}\doteq
c_{k-1}\pm 4b_{k-1}.  \tag{$\ddag $}
\end{equation}

Az $\left( \dag \right) $ és $\left( \ddag \right) $ miatt $%
a_{k},b_{k},c_{k}\in \mathbb{Z}$ valóban egész számok. Azt kell még
megmutatnunk, hogy $3$ nem osztja a $b_{k}$ egész számot. Ha $k=2,$ akkor $%
b_{2}=2\pm 2$ vagy $b_{2}=2\mp 0,$ tehát $b_{2}$ csak $0,2$ vagy $4$ lehet. N%
ézzük most a $k>2$ esetet! Négy esetet fogunk szétválasztani. \newline
1. Ha $w_{k}=\phi ^{\pm 1}\rho ^{\pm 1}w_{k-2}$ alakú, akkor $\left( \dag
\right) $ miatt $b_{k}=b_{k-1}\pm 2a_{k-1}$ és $\left( \ddag \right) $
szerint $b_{k}=b_{k-1}\pm 2\cdot 3a_{k-2}.$ Így az indukciós feltevés
szerint $3$ nem osztója $b_{k-1}$-nek ezért nem osztója $b_{k}$-nak sem.%
\newline
2. Ha $w_{k}=\rho ^{\pm 1}\phi ^{\pm 1}w_{k-2}$ alakban van felírva, akkor $%
\left( \ddag \right) $ miatt $b_{k}=b_{k-1}\mp 2c_{k-1}$ és $\left( \dag
\right) $ figyelembevételével $b_{k}=b_{k-1}\mp 2\cdot 3c_{k-2},$ ezért $3$
ebben az esetben sem osztója $b_{k}$-nak.\newline
3. Ha $w_{k}=\phi ^{\pm 1}\phi ^{\pm 1}w_{k-2}$ alakú, akkor $\left( \dag
\right) $-et kell kétszer alkalmaznunk, így 
\begin{eqnarray*}
b_{k} &=&b_{k-1}\pm 2a_{k-1}=b_{k-1}\pm 2\left( a_{k-2}\mp 4b_{k-2}\right)
=b_{k-1}\pm 2a_{k-2}-8b_{k-2} \\
&=&b_{k-1}+\left( b_{k-2}\pm 2a_{k-2}\right) -9b_{k-2}=2b_{k-1}-9b_{k-2}.
\end{eqnarray*}
Itt az utolsó lépésben $\left( \dag \right) $-et alkalmaztuk harmadszorra is 
$k$ helyett a $k-1$ számra. Mivel az indukciós feltevés szerint $3$ nem osztó%
ja $b_{k-1}$-nek, ezért $3$ nem osztója $b_{k}$-nak sem.\newline
4. Ha $w_{k}=\rho ^{\pm 1}\rho ^{\pm 1}w_{k-2}$ alakban van felírva, akkor
most $\left( \ddag \right) $-et kell kétszer alkalmaznunk, és úgy átalakí%
tani, hogy a $\left( \ddag \right) $ újbóli alkalmazásával az előzőhöz hasonl%
ó eset álljon elő. 
\begin{eqnarray*}
b_{k} &=&b_{k-1}\mp 2c_{k-1}=b_{k-1}\mp 2\left( c_{k-2}\pm 4b_{k-2}\right)
=b_{k-1}\mp 2c_{k-2}-8b_{k-2} \\
&=&b_{k-1}+\left( b_{k-2}\mp 2c_{k-2}\right) -9b_{k-2}=2b_{k-1}-9b_{k-2.}
\end{eqnarray*}
Beláttuk tehát, hogy $3$ nem osztója $b_{k}$-nak egyik esetben sem.
\end{proof}

Érdekes látni, hogy bizonyításunk egyedül itt lenne hibás, ha azt $\mathbb{R}%
^{2}$-re próbálnánk alkalmazni, hiszen a sík egységkörének forgatásai nyilvá%
nvaló módon kommutatív csoportot alkotnak, márpedig az triviális, hogy egy
kommutatív csoportnak csak az egy-egy elem generálta részcsoportjai
lehetnek szabad csoportok.

\begin{Allitas}[Hausdorff - 1914]
\index{Hausdorff-paradoxon} A $\mathbb{R}^{3}$ gömbfelszínének létezik olyan
megszámlálható részhalmaza, amelynek komplementere $SO\left( 3\right) $%
-paradox.
\end{Allitas}

\begin{proof}
Legyen $G$ az $SO\left( 3\right) $ két elem által generált szabad ré%
szcsoportja, amely (\ref{2szabad}) miatt paradox. Világos, hogy a négy jelbő%
l kirakható szavak halmaza legfeljebb megszámlálható számosságú és minden
forgatás, amely az identitástól különbözik éppen két pontot hagy érintetlenü%
l. Jelölje tehát $D$ az alábbi megszámlálható számosságú halmazt: 
\[
D\doteq \left\{ x\in S^{2}:\exists g\in G%
\text{ }g\left( x\right) =x\ \text{és}\ g\neq \id\right\} 
\]
Ha $g\in G$ esetén létezik $x\in S^{2}\smallsetminus D,$ melyre $g\left(
x\right) =x,$ akkor az csak úgy lehetséges, ha $g=\id$ is teljesül. Ha 
$x\in S^{2}$ esetén $g\left( x\right) \in D,$ akkor létezik $h\in G,h\neq 
\id,$ amelyre $h\left( g\left( x\right) \right) =g\left( x\right) ,$
ezért $x=g^{-1}hg\left( x\right) ,$ ahol $g^{-1}hg\neq \id.$
Eszerint $x\in D\Rightarrow g\left( x\right) \in D$ is fennáll minden $g\in G
$ esetén, emiatt $g\left( S^{2}\smallsetminus D\right) =S^{2}\smallsetminus D
$ is teljesül. A $G$ paradox csoport tehát $S^{2}\smallsetminus D$-n hat, ső%
t fixpontmentesen hat. Alkalmazható ezért az (\ref{hausdorff0}) állítás,
amely szerint az $S^{2}\smallsetminus D$ halmaz $G$-paradox, tehát $%
SO\left( 3\right) $-paradox is.
\end{proof}

\chapter{A Banach\texorpdfstring{\,--\,}{--}Tarski-paradoxon}

Hausdorff a fenti állítás segítségével megmutatta, hogy nincs az $S^{2}$ gö%
mbfelület egész hatványhalmazán értelmezett $\mu $ forgatás-invariáns vé%
gesen additív halmazfüggvény, melyre $\mu \left( S^{2}\right) =1.$ A bizonyít%
ása szerint egy ilyen $\mu $ minden megszámlálható halmazt $0$-ba visz, így
az állítás nyilvánvaló következménye a fenti Hausdorff-paradoxonnak.

Banach egyrészt megmutatta, hogy a Hausdorff-paradoxonból a megszámlálható
halmaz kihagyható, ezért rögtön látszik, hogy ilyen mérték valóban nincsen. M%
ásrészt észrevette, hogy könnyen áttérhetünk $B^{3}$ részhalmazainak vizsgá%
latára, sőt $\mathbb{R}^{3}$ két olyan részhalmaza, amely egyrészt nem túl nagy
(belefér egy gömbbe), másrészt nem túl kicsi (beleírható egy gömb) szintén
egymásba darabolható.

\section{A gyenge alak}

A Hausdorff-paradoxon megszámlálható halmazát az alábbi konkrét transzformá%
ciókkal fogjuk ignorálni:

\begin{Lemma}\label{btlemma}
    \begin{enumerate}
        \item 
Legyen $P\subseteq S^{1}$ egy megszámlálható halmaz é%
s $\tau _{\alpha }:\mathbb{R}^{2}\rightarrow \mathbb{R}^{2};\tau _{\alpha }\in
SO\left( 2\right) $ az $\alpha $ szögű forgatás. Ekkor megszámlálhatóan sok $%
\alpha $ kivételével minden $\alpha \in \left[ 0,2\pi \right) $ esetén 
\[
\left( \cup _{n=1}^{\infty }\tau _{\alpha }^{n}\left( P\right) \right) \cap
P=\varnothing .
\]
\item Legyen $P\subseteq S^{2}$ egy megszámlálható halmaz. Ekkor létezik $%
Q$ megszámlálható halmaz és $\omega :\mathbb{R}^{3}\rightarrow \mathbb{R}%
^{3};\omega \in SO\left( 3\right) $ forgatás, hogy 
\[
P\subseteq Q\subseteq S^{2}\quad \text{és}\quad \omega \left( Q\right)
=Q\smallsetminus P.
\]
\item Létezik $0\in N\subseteq B^{3}$ megszámlálható halmaz, és létezik $%
r:\mathbb{R}^{3}\rightarrow \mathbb{R}^{3};r\in G_{3}$ egybevágósági transzformáció%
, amelyre 
\[
r\left( N\right) =N\smallsetminus \left\{ 0\right\} .
\]
Alkalmazva tehát az (\ref{tolo}) állítást azt kapjuk, hogy amennyiben $P$
megszámlálható, akkor 
\[
S^{2}\sim _{SO\left( 3\right) }\left( S^{2}\smallsetminus P\right) \quad 
\text{valamint\quad }B^{3}\sim _{G_{3}}\left( B^{3}\smallsetminus \left\{
0\right\} \right) .
\]
    \end{enumerate}
\end{Lemma}

\begin{proof}
\indent
$\left( 1\right) $ Valamely $\alpha \in \left[ 0,2\pi \right) $ mellett
legyen $z_{\alpha }\doteq e^{i\alpha }.$ Ha erre a az $\alpha $-ra és valamely 
$n\in \mathbb{N}$-re $\tau _{\alpha }^{n}\left( P\right) \cap P\neq \varnothing 
$ azt jelenti, hogy létezik $x,y\in P,$ amelyre $z_{\alpha }^{n}x=y,$ azaz $%
z_{\alpha }=\sqrt[n]{\frac{y}{x}}.$ Mivel egy komplex számnak csak $n$ darab 
$n$-edik gyöke van, ez azt jelenti, hogy adott $\left( x,y,n\right) \in
P\times P\times \mathbb{N}$-hez, csak véges sok $\alpha $ létezik, amelyre $%
\tau _{\alpha }^{n}\left( P\right) \cap P\neq \varnothing .$ Megismételve
ezt $P\times P\times \mathbb{N}$ minden elemére, azt kapjuk, hogy csak megszámlálhatóan 
sok olyan $\alpha \in \left[ 0,2\pi \right) $ szög van, amelyre $%
\left( \cup _{n=1}^{\infty }\tau _{\alpha }^{n}\left( P\right) \right) \cap
P\neq \varnothing .$

$\left( 2\right) $ Legyen most $P\subseteq S^{2}\subseteq \mathbb{R}^{3}$ egy
megszámlálható halmaz. Válasszunk egy olyan $L$ egyenest, amely $S^{2}$-et $P
$-től különböző pontokban metszi. Mivel $P$ csak megszámlálható számosságú,
ilyen egyenes biztosan van. Legyen $\omega _{\alpha }:\mathbb{R}^{3}\rightarrow 
\mathbb{R}^{3}$ ezen $L$ egyenes mint tengely körüli $\alpha $ szögű forgatás.
Világos, hogy csak megszámlálható olyan sík van, amely $L$-re merőleges és
tartalmaz $P$-beli pontot. Erre a legfeljebb megszámlálhatóan sok sí%
kmetszetre alkalmazva a lemma első felében már igazolt állítást, azt kapjuk,
hogy megszámlálhatóan sok $\alpha $ szög kivételével tetszőleges $\alpha $%
-ra 
\[
\left( \cup _{n=1}^{\infty }\omega _{\alpha }^{n}\left( P\right) \right)
\cap P=\varnothing .
\]
Egy ilyen $\alpha $ mellett tekintsük a $Q\doteq P\cup \left( \cup
_{n=1}^{\infty }\omega _{\alpha }^{n}\left( P\right) \right) $ diszjunkt
egyesítést. Világos, hogy $P\subseteq Q\subseteq S^{2};$ $Q$ megszámlálható
számosságú; és $\omega \left( Q\right) =\cup _{n=1}^{\infty }\omega _{\alpha
}^{n}\left( P\right) =Q\smallsetminus P.$

$\left( 3\right) $ Jelölje $u\doteq \left( 1,0,0\right) \in S^{2}$. A lemma m%
ásodik állítását alkalmazva legyen $u\in Q\subseteq S^{2}$ megszámlálható
halmaz, valamint $\rho $ az forgatás, amelyre 
\[
\rho \left( Q\right) =Q\smallsetminus \left\{ u\right\} .
\]
Legyen $N\doteq \frac{1}{2}Q-\frac{1}{2}u.$ Világos, hogy $N\subseteq B$,
amelyre $0\in N$. Definiáljuk az alábbi $r$ egybevágósági transzformációt: 
\[
r\left( x\right) =\rho \left( x+\frac{1}{2}u\right) -\frac{1}{2}u.
\]
Ha $x\in N$, akkor $x=\frac{1}{2}q-\frac{1}{2}u$ valamely $q\in Q$-val, így $%
r\left( x\right) =\rho \left( \frac{1}{2}q-\frac{1}{2}u+\frac{1}{2}u\right) -%
\frac{1}{2}u=\frac{1}{2}\rho \left( q\right) -\frac{1}{2}u,$ amiből 
\[
r\left( N\right) =\frac{1}{2}\rho \left( Q\right) -\frac{1}{2}u=\frac{1}{2}%
\left( Q\smallsetminus \left\{ u\right\} \right) -\frac{1}{2}%
u=N\smallsetminus \left\{ 0\right\} .
\]
\end{proof}

\begin{Allitas}[Banach\,--\,Tarski-paradoxon gyenge alakja]
\index{Banach\,--\,Tarski-paradoxon gyenge alak} Az $S^{2}$ gömbfelület $%
SO\left( 3\right) $-paradox, valamint a $B^{3}$ gömb $G_{3}$-paradox.
\end{Allitas}

\begin{proof}
A Hausdorff-paradoxon szerint $S^{2}\smallsetminus P$ egy $SO\left(
3\right) $-paradox halmaza a gömbfelszínnek. Az előző lemma miatt viszont$\
\left( S^{2}\smallsetminus P\right) \sim _{SO\left( 3\right) }S^{2}$, ezért
a (\ref{partranz}) állítást alkalmazva azt kapjuk, hogy $S^{2}$ valóban $%
SO\left( 3\right) $-paradox halmaz.

A gömbfelület egy részhalmazához hozzárendelve ugyanezen gömbfelületrészhez
tartozó gömbcikket, könnyen látható, hogy $B^{3}\smallsetminus \left\{
0\right\} $ pedig $SO\left( 3\right) $-paradox részhalmaza, ezért $G_{3}$%
-paradox részhalmaza is a $B^{3}$ gömbnek. Az előző lemma miatt viszont $%
B^{3}\smallsetminus \left\{ 0\right\} \sim _{G_{3}}B^{3},$ ezért a fentihez
hasonlóan a (\ref{partranz}) állítást alkalmazva kapjuk, hogy $B^{3}$ való%
ban egy $G_{3}$-paradox halmaz.
\end{proof}

\section{Az erős alak}

Az erős alak igazolásának egyetlen trükkje, hogy a számosságok elméletéből jó%
l ismert Schröder\,--\,Bernstein-ekvivalenciatételt kell élesítenünk bijekciók
helyett távolságtartó bijekciókra. Ehhez először bevezetünk egy relációt,
amelynek szimmetrikus része az ekvidekompozábilitás relációja lesz.

\begin{Definicio}[ekvidekompozábilis beágyazás]
\index{ekvidekompozábilis beágyazás} Tegyük fel, hogy $G$ egy az $X$
halmazon ható csoport, $A,B\subseteq X.$ Azt mondjuk, hogy $A$ ekvidekompozá%
bilis része $B$-nek, ha létezik $B^{\prime }\subseteq B,$ amelyre $A\sim
_{G_{3}}B^{\prime }.$ Ezt $A\preceq _{G}B$ módon jelöljük.
\end{Definicio}

\begin{Allitas}[Banach\,--\,Schröder\,--\,Bernstein]
\index{Banach\,--\,Schröder\,--,Bernstein-ekvivalenciatétel} Legyen $G$ egy az $X$
halmazon ható csoport és $A,B\subseteq X.$ Ekkor $B\preceq _{G}A$ és $%
A\preceq _{G}B$ esetén $A\sim _{G}B$ is fennáll.
\end{Allitas}

\begin{proof}
Legyen $A=\cup _{i=1}^{n}A_{i}$ diszjunkt felbontás és $g_{i}\in G,$
amelyekre $\cup _{i=1}^{n}%
\hat{g}_{i}\left( A_{i}\right) =B^{\prime }\subseteq B$ egy partíciója $%
B^{\prime }$-nek. Hasonlóan, $B=\cup _{j=1}^{m}B_{j}$ partíció és $f_{j}\in
G,$ amelyekre $\cup _{j=1}^{m}\hat{f}_{j}\left( B_{j}\right) =A^{\prime
}\subseteq A$ is partíciója $A^{\prime }$-nek. Legyenek a $g:A\rightarrow
B^{\prime },g\left( x\right) \doteq \hat{g}_{i}\left( x\right) ,$ $x\in
A_{i} $ és az $f:B\rightarrow A^{\prime },f\left( x\right) \doteq \hat{f}%
_{j}\left( x\right) ,x\in B_{j}$ módon definiált bijekciók. Tekintsük az 
\[
\Phi :\mathcal{P}\left( A\right) \rightarrow \mathcal{P}\left( A\right)
;\qquad \Phi \left( D\right) \doteq A\smallsetminus f\left( B\smallsetminus
g\left( D\right) \right) 
\]
halmaz-halmaz leképezést. Világos, hogy $D_{1}\subseteq D_{2}$ esetén $\Phi
\left( D_{1}\right) \subseteq \Phi \left( D_{2}\right) $ teljesül.
Megmutatjuk, hogy létezik $\Phi $-nek fixpontja. Ehhez tekintsük a 
\[
\mathcal{H}\doteq \left\{ C\subseteq A:C\subseteq \Phi \left( C\right)
\right\} 
\]
halmazrendszert. Mivel $\varnothing \in \mathcal{H,}$ ezért $\mathcal{H}$
nem üres, így legyen $D\doteq \cup \mathcal{H}$ az összes $\mathcal{H}$-beli
halmaz egyesítése. Persze $C\in \mathcal{H}$ mellett $C\subseteq D,$ ezért $%
C\subseteq \Phi \left( C\right) \subseteq \Phi \left( D\right) ,$ amiből $%
D\subseteq \Phi \left( D\right) $ következik. Így tehát $\Phi \left(
D\right) \subseteq \Phi \left( \Phi \left( D\right) \right) ,$ ezért $\Phi
\left( D\right) \in \mathcal{H,}$ ergo $\Phi \left( D\right) \subseteq D$ is
teljesül. Találtunk tehát $D\subseteq X$ halmazt, amelyre 
\begin{equation}
A\smallsetminus D=f\left( B\smallsetminus g\left( D\right) \right) 
\tag{$\dag $}
\end{equation}
teljesül. Legyen $i=1,\ldots ,n$ és $j=1,\ldots ,m$ esetén 
\[
E_{i}\doteq A_{i}\cap D;\quad h_{i}\doteq g_{i};\quad E_{n+j}\doteq \hat{f}%
_{j}\left( B_{j}\smallsetminus g\left( D\right) \right) ;\quad h_{n+j}\doteq
f_{j}^{-1}. 
\]
Világos, hogy $D=\cup _{i=1}^{n}E_{i}$ partíciója $D$-nek; $B\smallsetminus
g\left( D\right) =\cup _{j=}^{m}\left( B_{j}\smallsetminus g\left( D\right)
\right) $ partíciója $B\smallsetminus g\left( D\right) $-nek, ezért $\left(
\dag \right) $ miatt $A\smallsetminus D=\cup _{j=}^{m}\hat{f}_{j}\left(
B_{j}\smallsetminus g\left( D\right) \right) =\cup _{j=1}^{m}E_{n+j}$ partíciója 
$A\smallsetminus D$-nek. Együtt tehát 
\[
A=\cup _{k=1}^{n+m}E_{k} 
\]
partíciója $A$-nak. \newline
Mivel $\cup _{i=1}^{n}\hat{g}_{i}\left( A_{i}\right) =B^{\prime }\subseteq B$
partíciója $B^{\prime }$-nek, ezért $\cup _{k=1}^{n}\hat{h}_{k}\left(
E_{k}\right) =\cup _{i=1}^{n}\hat{g}_{i}\left( A_{i}\cap D\right) =g\left(
D\right) \subseteq B^{\prime }$ partíciója $g\left( D\right) $-nek. Világos,
hogy $\hat{h}_{n+j}\left( E_{n+j}\right) =B_{j}\smallsetminus g\left(
D\right) ,$ ezért $\cup _{k=n+1}^{n+m}\hat{h}_{k}\left( E_{k}\right) $ pedig
partíciója $B\smallsetminus g\left( D\right) $-nek. Azt kaptuk tehát, hogy 
\[
\cup _{k=1}^{n+m}\hat{h}_{k}\left( E_{k}\right) 
\]
partíciója $B$-nek. Ezt kellett belátni.
\end{proof}

Első következményként egyszerűsíthetjük a paradox halmaz definícióját.

\begin{Allitas}[paradox halmaz felbontásáról]
\label{paradox}Legyen $G$ egy az $X$ halmazon ható csoport. Legyen $%
C\subseteq A$ olyan halmazok, amelyekre $C\sim _{G}A.$ Ekkor tetszőleges $%
C\subseteq Q\subseteq A$ esetén $Q\sim _{G}A$ is fennáll. \newline
Az $A\subseteq X$ halmaz pontosan akkor $G$-paradox, ha $A$ előáll $A=B\cup
C$ diszjunkt egyesítés alakban, ahol $B$ és $C$ is $G$-ekvidekompozábilisek 

$A$-val.
\end{Allitas}

\begin{proof}
Világos, hogy a tartalmazások miatt $C\preceq _{G}Q\preceq _{G}A.$ No de $%
A\preceq _{G}C$ is fennáll, ezért a Banach\,--\,Schröder\,--\,Bersnstein-]ekvivalenciat%
étel miatt $Q\sim _{G}A$ is teljesül. A második állítás már az első következm%
énye.
\end{proof}

\begin{Allitas}[Diszjunkt duplázás]
Legyen $A\subseteq \mathbb{R}^{3}$ egy zárt gömb és $A_{1}$ ennek olyan
eltoltja, hogy $A\cap A_{1}=\varnothing .$ Ekkor $A\sim _{G_{3}}A\cup A_{1}.$
\end{Allitas}

\begin{proof}
A gyenge alak és az előző állítás szerint az $A$ gömbnek van olyan $A=B\cup
C $ diszjunkt előállítása, ahol $B$ és a $C$ részhalmaz is $G_{3}$%
-ekvidekompozábilis $A$-val. Ekkor persze $C$ is $G_{3}$-ekvidekompozábilis 
$A^{\prime }$-vel. Mivel $A^{\prime }$ és $A$ diszjunktak, ezért $%
A=B\cup C\sim _{G_{3}}A\cup A^{\prime }.$
\end{proof}

\begin{Allitas}[$n$-szerezés]
Legyen $A\subseteq \mathbb{R}^{3}$ egy zárt gömb, és $A_{1},\ldots ,A_{n}$
ennek véges számú eltoltja. Ekkor 
\[
A\sim _{G_{3}}\cup _{j=1}^{n}A_{j} 
\]
\end{Allitas}

\begin{proof}
Persze $n$ szerinti teljes indukció, és az $n=1$ eset csak azt állítja, hogy 
$A\sim _{G_{3}}A$. Tegyük fel, hogy, hogy $A\sim \cup _{j=1}^{n-1}A_{j}$ és l%
ássuk be $n$-re. Világos, hogy $A_{n}\smallsetminus \cup _{j=1}^{n-1}A_{j}$
az $A_{n}$ részhalmazaként ekvidekompozábilis $A$ egy részhalmazával, tehát $%
A_{n}\smallsetminus \cup _{j=1}^{n-1}A_{j}\preceq _{G_{3}}B$, ahol $B$ az $A$%
-nak olyan eltolja, melyre $A\cap B=\varnothing .$ Ezért 
\[
\cup _{j=1}^{n}A_{j}=\left( \cup _{j=1}^{n-1}A_{j}\right) \cup \left(
A_{n}\smallsetminus \cup _{j=1}^{n-1}A_{j}\right) \preceq _{G_{3}}A\cup
B\sim _{G_{3}}A. 
\]
No de, $A\preceq \cup _{j=1}^{n}A_{j}$ triviális, ezért az ekvivalenciatétel
miatt $A\sim \cup _{j=1}^{n}A_{j}$ is valóban teljesül.
\end{proof}

\begin{Allitas}[Banach\,--\,Tarski-paradoxon erős alakja]
\index{Banach\,--\,Tarski-paradoxon erős alak} Legyen $Q$ és $T$ az $\mathbb{R}^{3}$
egy-egy korlátos részhalmaza, melyeknek belseje sem üres. Ekkor $P\sim
_{G_{3}}T$, azaz a $P$ és a $T$ halmazok $G_{3}$-ekvidekompozábilisek egymással.
\end{Allitas}

\begin{proof}
Legyen $A$ egy olyan zárt gömb, melyre $A\subseteq Q$ és melynek bizonyos $B$
eltoltjára $B\subseteq T.$ Mivel $Q$ teljesen korlátos, ezért vannak $A$-nak
olyan $A_{1},\ldots ,A_{n}$ eltoltjai, amelyekre $Q\subseteq \cup
_{j=1}^{n}A_{j}.$ Világos, hogy a tartalmazások és az előző tétel szerint 
\[
A\preceq _{G_{3}}Q\preceq _{G_{3}}\cup _{j=1}^{n}A_{j}\sim _{G_{3}}A, 
\]
így újra az ekvivalencia tételt használva azt kapjuk, hogy $A\sim _{G_{3}}Q.$
Ugyanezt kapjuk $A$-nak a bizonyítás elején definiált $B$ eltoltjára és $T$%
-re, azaz $B\sim _{G_{3}}T.$ Összevetve: $Q\sim _{G_{3}}A\sim _{G_{3}}B\sim
_{G_{3}}T$ valóban teljesül.
\end{proof}

\section*{Köszönetek}

Köszönettel tartozom Dancs Istvánnak, aki az első elemi bizonyításra
\cite{Stromberg-1979} felhívta figyelmem. 
Szintén köszönettel tartozom
Petrus Potgieternek barátomnak \cite{Schechter-2000} felkutatásáért és a halmazelméleti diszkussziókért.

\bibliographystyle{plain}
\bibliography{btarxiv.bib}

\bigskip
\hfill{\ 
\begin{tabular}{l}
{Delians: ``How can we be rid of the plague?''} \\ 
{Delphic Oracle: ``Construct a cubic altar of the double } \\ 
\qquad {the size of the existing one.''} \\ 
{Banach and Tarski: ``Can we use the Axiom of Choice?''}
\end{tabular}
}

\hfill{Stan Wagon (1985) \cite{Wagon-1985}}

\bigskip
\hfill{,,Ők így válaszoltak: Nincs itt egyebünk, csak öt kenyerünk és két halunk.''}

\hfill{Mt 14,17}

\end{document}